\documentclass[11pt,a4paper]{amsart}

\pdfoutput=1
\usepackage[english]{babel}
\usepackage{graphicx}
\usepackage{framed}
\usepackage[normalem]{ulem}
\usepackage{amsmath}
\usepackage{amsthm}
\usepackage{amssymb}
\usepackage{amsfonts}
\usepackage[foot]{amsaddr}
\usepackage{mathrsfs,amsmath}
\usepackage{enumerate}
\usepackage[utf8]{inputenc}
\usepackage[top=1 in,bottom=1in, left=1 in, right=1 in]{geometry}
\usepackage{mathtools}
\usepackage{tikz-cd}
\usepackage[font=small,labelfont=bf]{caption} 
\usepackage[subrefformat=parens]{subcaption}
\usepackage{cite}
\usepackage{abstract}

\newcommand{\area}{{\rm area}}

\newcommand{\rk}{{\rm rk}}

\newcommand{\R}{\mathbb{R}}
\newcommand{\N}{\mathbb{N}}

\newcommand{\Z}{\mathbb{Z}}

\newcommand{\<}{\langle}

\renewcommand{\>}{\rangle}
\newcommand{\eps}{\varepsilon}

\newtheorem{thm}{Theorem}[section]

\newtheorem{lem}[thm]{Lemma}

\theoremstyle{definition}
\newtheorem{defn}[thm]{Definition}

\theoremstyle{remark}

\newtheorem{rmk}[thm]{Remark}

\newcommand{\eproof}{\hfill\qed}

\newcommand{\bproofof}[1]{\noindent{\textit{Proof of #1. }}}

\title{$C^0$-limits of Legendrian Submanifolds}
\author{Lukas Nakamura}
\email[A1]{nakamura@math.lmu.de}
\date{}

\begin{document}

\maketitle

\begin{abstract}
Laudenbach and Sikorav proved that closed, half-dimensional non-Lagrangian submanifolds of symplectic manifolds are immediately displaceable as long as there is no topological obstruction. From this they deduced that under certain assumptions the $C^0$-limit of a sequence of Lagrangian submanifolds is again Lagrangian, provided that the limit is smooth. 

In this note we extend Laudenbach and Sikorav's ideas to contact manifolds. We prove correspondingly that certain non-Legendrian submanifolds of contact manifolds can be displaced immediately without creating short Reeb chords as long as there is no topological obstruction. From this it will follow that under certain assumptions the $C^0$-limit of a sequence of Legendrian submanifolds with uniformly bounded Reeb chords is again Legendrian, provided that the limit is smooth. 
\end{abstract}

\section{Introduction}

The Lagrangian Arnold conjecture \cite{arn65} implies that a Lagrangian submanifold $L$ of a symplectic manifold $M$ always intersects its image under a Hamiltonian diffeomorphism. Furthermore, the number of intersection points should be bounded from below by the Betti number of $L$ if the intersection is transverse and by the cup-length of $L$ in the general case. For $C^1$-small Hamiltonian diffeomorphisms in a cotangent bundle the Arnold conjecture follows easily from Morse theory. Gromov proved in \cite{gro85} with the use of pseudo-holomorphic curves that it is impossible to displace a weakly exact Lagrangian in a geometrically bounded symplectic manifold by a Hamiltonian diffeomorphism. Of course, this cannot hold for arbitrary Lagrangians in arbitrary symplectic manifolds as the example of an embedded circle in $\R^2$ with its standard symplectic structure shows. But Polterovich \cite{pol93} showed under the assumptions that $L$ is rational and that $M$ is geometrically bounded that $L$ will always intersect its image under a Hamiltonian diffeomorphism $\psi$ as long as $\psi$ is sufficiently small in the Hofer norm. Floer \cite{flo88} introduced a homology theory for Lagrangian intersections in order to prove the Arnold conjecture in the case that $M$ is compact and $\pi_2(M,L) = 0$. Chekanov \cite{che98} used Floer's ideas to prove that the Arnold conjecture holds for all closed Lagrangians in geometrically bounded symplectic manifolds as long as the Hamiltonian diffeomorphism is sufficiently small in the Hofer norm.
 
Now, the question arises whether non-Lagrangian submanifolds can be rigid as well. To this end, Laudenbach and Sikorav proved in \cite{ls94} that half-dimensional closed non-Lagrangian submanifolds of symplectic manifolds are infinitesimally displaceable as long there is no topological obstruction. Here, infinitesimally displaceable means that there is a Hamiltonian vector field nowhere tangent to that submanifold.

Similarly to the symplectic case, there are also results about the rigidity of Legendrian submanifolds $L$ in a contact manifold $M$. For example,  Rizell and Sullivan (\!\!\cite{rs16}, \cite{rs18}) proved that if the contact Hamiltonian $H$ generating a contactomorphism $\phi^H$ is ``sufficiently small", then there are short (compared to $H$) Reeb chords between $L$ and $\phi(L)$.

In this work, we extend Laudenbach and Sikorav's ideas to contact manifolds. We prove that under certain assumptions for a given $n$-dimensional non-Legendrian submanifold $L$ (where $\dim(M) = 2n+1$) there exists a contact vector field that is nowhere contained in the sum of the tangent space of $L$ and the span of the Reeb vector field along $L$.

Laudenbach and Sikorav \cite{ls94} noted that if a sequence $\{L_n\}_{n \in \N}$ of closed Lagrangian submanifolds of a geometrically bounded\footnote{They consider the cases $M = \R^{2n}$ and $\pi_2(M,L) = 0$ but their proof easily extends to general geometrically bounded symplectic manifolds, cf. Theorem \ref{c0-limit lagrange general} below.} symplectic manifold $C^0$-converges to an embedded submanifold $L$, then the displacement energies of the $L_n$ have to be uniformly bounded away from zero. But if $L$ has vanishing displacement energy, then the sequence of the displacement energies of the $L_i$ has to go to zero. From this they concluded that the limit has to be Lagrangian as well. 

In a similar way, it will follow that the limit of a sequence of closed Legendrian submanifolds with uniformly bounded Reeb chords is again Legendrian (Theorem \ref{c0-limit legendre}).\\

\textit{Acknowledgements}: This work was carried out as part of the Master's program ``Theoretical and Mathematical Physics" at the Ludwig-Maximilans-University Munich, and it summarizes the results of my Master's thesis. I would first like to thank Thomas Vogel for supervising this work and for his numerous helpful remarks about this note. He always found the time to answer all of my questions. Furthermore, I am grateful to Yang Huang for many interesting and stimulating discussions. Also, I would like to thank Georgios Dimitroglou Rizell for explaining to me some of the results of his joint work with M. Sullivan.

\hfill \break
\section{Displacing non-Legendrian submanifolds}\label{sec:main results}

As mentioned in the introduction, closed Lagrangian submanifolds of many symplectic manifolds are rigid. Let us describe the following rather weak rigidity property. Let $(M^{2n},\omega)$ be a symplectic manifold and $L \subseteq M$ a closed Lagrangian submanifold. The restriction of any function $H:M \to \R$ to $L$ has a critical point $x \in L$ because $L$ is closed, i.e. $dH(x)|_{T_xL} = 0$. For the Hamiltonian vector field $X_H$ associated to $H$, defined by $i_{X_H} \omega = - dH$, this  implies that $X_H(x) \in T_xL^{\perp_\omega} = T_xL$ since $L$ is Lagrangian. In other words, there exists no Hamiltonian vector field on $M$ that is nowhere tangent to $L$.  

Now let $L^n$ be a closed non-Lagrangian submanifold of $M$ and we ask whether there exists a Hamiltonian vector field nowhere tangent to $L$. Of course, there might not exist any vector field that is nowhere tangent to $L$ as the self-intersection number of $L$ might be non-zero. But under the additional assumption that there is no such topological obstruction, Laudenbach and Sikorav proved the affirmative answer.

\begin{thm}\label{main thm symp}\hspace{-1mm}\normalfont{\cite{ls94}}\,
Let ($M^{2n},\omega$) be a symplectic manifold and $L$ a closed, connected submanifold of dimension n such that

(i) $L$ is $non$-$Lagrangian$, i.e. there exists a point $x \in L$ such that $T_xL$ is not a Lagrangian subspace of $T_xM$,

(ii) the normal bundle $\nu$ of $L \subseteq M$ has a nowhere vanishing section. 

Then there exists a Hamiltonian vector field on $M$ that is nowhere tangent to $L$.\\
\end{thm}

\begin{rmk}\label{rmk:generalizations of symp non-rig}
Clearly, the generalization of Theorem \ref{main thm symp} to non-coisotropic submanifolds fails in general as such manifolds may contain closed Lagrangian submanifolds. However, Gürel \cite{gue08} noted that Theorem \ref{main thm symp} extends to nowhere coisotropic manifolds. Also, one can prove that even the parametric and a relative version of the h-principle for Hamiltonian vector fields that are nowhere tangent to $L$ holds.\\
\end{rmk}

Analogously to the Lagrangian case, Legendrians obey the following rigidity result. Let $(M, \xi = \ker \alpha)$ be a cooriented contact manifold and $L \subseteq M$ a closed Legendrian submanifold. Let $H: M \to \R$ be an arbitrary function. Then $H|_L$ has a critical point $x \in L$. From $dH(x)|_{T_xL} = 0$ it follows that $X_H(x) \in T_xL^{\perp_{d \alpha}} \oplus \< R_\alpha(x) \> = T_xL \oplus \< R_\alpha(x) \>$. Here, $X_H$ denotes the contact vector field associated to $H$ that is defined by
\begin{equation}
i_{X_H} d \alpha|_{\xi} = -dH|_\xi \quad \text{ and } \quad \alpha(X_H) = H,
\end{equation}
$(\cdot)^{\perp_{d \alpha}}$ denotes the $d \alpha|_\xi$ complement in $\xi$, and $R_\alpha$ denotes the Reeb vector field on $M$. This means that for a closed Legendrian submanifold there exists no contact vector field that is nowhere contained in $TL \oplus R_\alpha$.

We now also consider non-Legendrian submanifolds. Below, we will apply the proof of Theorem \ref{main thm symp} in \cite{ls94} in order to show that, as in the symplectic case, there exist contact vector fields nowhere tangent to $TL \oplus R_\alpha$ as long as there is no topological obstruction, at least for a generic non-Legendrian submanifold. Note that the flow of such a contact vector field displaces the non-Legendrian submanifold $L$ in such a way that there are no short Reeb chords between $L$ and its image under the flow. 

\begin{thm}\label{main thm contact}
Let ($M^{2n+1}, \xi = \mathrm{ker}\,\alpha$) be a cooriented contact manifold. Denote its Reeb vector field by $R_\alpha$. Let $L \subseteq M$ be a closed, connected submanifold of dimension n such that

(i) $R_\alpha(x) \notin T_x L$ for all $x \in L$,

(ii) $L$ is $non$-$Legendrian$, i.e. there exists a point $x \in L$ with $T_xL \not\subseteq \xi_x$,

(iii) there exists a nowhere vanishing section of the normal bundle of the subvector bundle $TL \oplus \< R_\alpha|_L \> \subseteq TM|_L$.

 Then there exists a contact vector field $X$ such that $X(x) \notin T_xL \,\oplus \< R_\alpha (x) \>$ for all $x \in L$.\\
\end{thm}

\begin{rmk}\label{rmk:conditions in cont non-rig are necessary}
For a generic $n$-dimensional submanifold $L \subseteq (M,\ker \alpha)$, $R_\alpha$ will be nowhere tangent to $L$. Hence, Theorem \ref{main thm contact} describes the generic case. With basically the same proof one can show that a similar statement also holds if we require that the Reeb vector field is everywhere tangent to $L$.\\
\end{rmk}

\begin{rmk}
Similarly to Gürel's result \cite{gue08} that was mentioned in Remark \ref{rmk:generalizations of symp non-rig}, Theorem \ref{main thm contact} also holds for submanifolds that have a dimension different from $n$ if one requires that $\left( \pi T_xL \right)^{\perp_{d \alpha}} \not\subseteq T_xL$ holds for all $x \in L$. Here, $\pi:TM = \xi \oplus \< R_\alpha \> \to \xi$ denotes the projection onto the first factor. Also, the relative and a parametric h-principle hold in the setting of Theorem \ref{main thm contact} and in this case.\\
\end{rmk}

Laudenbach and Sikorav deduced Theorem \ref{main thm symp} from the following more general statement.

\begin{thm}\label{tech thm}\hspace{-1mm}\normalfont{\cite{ls94}}\,
Let $M$ be a manifold, $L$ a closed connected submanifold, and $E$ a subbundle of $TM|_L$ with $\rk(E) = \dim(L)$ such that

(i) $E\neq TL$, i.e. there exists a point $x \in L$ with $E_x\neq T_xL$,

(ii) there exists a nowhere vanishing section of $E$.\\
Then there exists a function $H$ on $M$ such that $dH|_{E_x}$ is non-zero for all $x \in L$.\\
\end{thm} 

\bproofof{Theorem~\ref{main thm contact}}
We will show how Theorem \ref{tech thm} implies Theorem \ref{main thm contact}. 

The tangent bundle $TM$ of $M$ splits as $TM = \xi \oplus \< R_\alpha \>$. As above, let $\pi$ denote the projection onto the first factor. In order to apply Theorem \ref{tech thm}, we define the vector bundle 
\begin{equation}
E \coloneqq (\pi TL)^{\perp_{d \alpha}}
\end{equation}
on $L$. Since $R_\alpha$ is nowhere tangent to $L$, this indeed defines a vector bundle with $\rk(E) = \dim(L)$. Because $E \subseteq \xi$, it follows that $E = TL$ if and only if $L$ is Legendrian. Thus, condition $(ii)$ in Theorem \ref{main thm contact} is precisely condition $(i)$ in Theorem \ref{tech thm}.

It is convenient to consider a complex structure $J: \xi \to \xi$ on the contact distribution such that
\begin{equation}
g_J(v,w) \coloneqq d \alpha (v, Jw), \quad v,w \in \xi_x, x \in L,
\end{equation}
defines a metric on $\xi$. Such complex structures exist because $d \alpha|_\xi$ defines a symplectic structure on $\xi$ (cf. \cite{ms17}, Proposition 2.6.4). Then we can extend $g_J$ to a metric on $M$ in such a way that the Reeb vector field $R_\alpha$ is orthogonal to $\xi$. 

By assumption, there exists a vector field $X$ that is orthogonal to $TL \oplus \< R_\alpha \>$ at every point of $L$. Especially, $X$ is tangent to $\xi$ along $L$. Now it is easy to check that $J X$ defines a nowhere vanishing section of $E$. 

Therefore, Theorem \ref{tech thm} implies that there exists a function $H: M \to \R$ such that $dH|_{E_x}$ is non-zero for all $x \in L$. For any $x \in L$, we have that 
\begin{equation}
0 = dH|_{E_x} = dH|_{(\pi T_xL)^{\perp_{d \alpha}}} \quad \Leftrightarrow \quad X_H(x) \in T_xL \oplus \< R_\alpha(x) \>.
\end{equation}
Hence, Theorem \ref{main thm contact} follows. 
\eproof

\hfill \break
\section{$C^0$-limits of Legendrian submanifolds}\label{sec:applications}

Let $(M,\omega)$ be a symplectic manifold. Eliashberg \cite{eli87} proved that the group of symplectomorphisms of $M$ is $C^0$-closed as a subset of the group of diffeomorphisms of $M$. This theorem can be stated equivalently in terms of graphs of diffeomorphisms of $M$. For this, recall that a diffeomorphism of $M$ is a symplectomorphism if and only if its graph in $(M \times M, pr_1^* \omega - pr_2^* \omega)$ is Lagrangian. Then Eliashberg's result states that the $C^0$-limit of a sequence of smooth, Lagrangian graphs in $M \times M$ is again Lagrangian, provided that it is a smooth graph.

Now one can also consider the closure of the symplectomorphism group of $M$ inside the group of homeomorphisms of $M$. A homeomorphism that is a $C^0$-limit of symplectomorphisms is called a $C^0$-symplectomorphism. Humilière, Leclercq and Seyfaddini \cite{hls15} generalized Elishberg's Theorem: If a $C^0$-symplectomorphism maps a coisotropic submanifold to a smooth manifold, then the image will be coisotropic as well. 

In these statements it is assumed that the $C^0$-limits of the Lagrangian (or coisotropic) submanifolds are induced by $C^0$-limits of symplectomorphisms. But Laudenbach and Sikorav showed that this assumption is not necessary in general.

\begin{thm}\label{c0-limit lagrange}\hspace{-1mm}\normalfont{\cite{ls94}}\,
Let $(M^{2n}, \omega)$ be a symplectic manifold and $L^n$ a closed manifold. Let $f_i: L \to M$ be a sequence of Lagrangian embeddings of $L$ into $M$ that $C^0$-converges to an embedding $f:L \to M$. If 

(i) $(M, \omega)$ is geometrically bounded and $\pi_2(M, f(L)) = 0$, or

(ii) $(M, \omega) = (\R^{2n}, \omega_0)$,\\
then $f$ is a Lagrangian embedding.\\
\end{thm}

Recall the following definition.

\begin{defn}\label{geombdd} (cf. \cite{gro85}, \cite{al94}) Let $(M,\omega)$ be a symplectic manifold. It is \textit{geometrically bounded} if there exists an almost complex structure $J$ such that $g_J(\cdot,\cdot) \coloneqq \omega(\cdot,J \cdot)$ defines a complete Riemannian metric for which there exists an upper bound on the sectional curvature and a positive lower bound on the injectivity radius of $(M,g_J)$.\\
\end{defn}

We will show that Theorem \ref{c0-limit lagrange} even holds if we replace the conditions $(i)$ and $(ii)$ by the more general condition that $(M,\omega)$ is geometrically bounded.

\begin{thm}\label{c0-limit lagrange general}
Let $(M^{2n}, \omega)$ be a geometrically bounded symplectic manifold and $L^n$ a closed manifold. Let $f_i: L \to M$ be a sequence of Lagrangian embeddings of $L$ into $M$ that $C^0$-converges to an embedding $f:L \to M$. Then $f$ is a Lagrangian embedding.\\
\end{thm}

Now consider a cooriented contact manifold $(M, \ker \alpha)$. Correspondingly to Eliashberg's result, Müller and Spaeth \cite{ms14} showed that the group of contactomorphism of $M$ is $C^0$-closed as a subset of the group of diffeomorphisms. Again, we also obtain a statement about the graphs of contactomorphisms as follows. Consider the projections $pr_1, pr_2: E \coloneqq M \times M \times \R \to M$ onto the first and second factor, respectively. A section of $pr_1: (E, e^z pr_1^*\alpha - pr_2^*\alpha) \to M$ is Legendrian if and only if it is of the form $x \mapsto (x, \psi(x), g(x))$ for some contactomorphism $\psi:M \to M$. Here, $z$ denotes the coordinate on $\R$ and $g$ is the conformal factor of $\psi$ defined by $\psi^*\alpha = e^g \alpha$. If we now apply Müller and Spaeth's Theorem to a sequence of contactomorphisms for which their respective conformal factors converge uniformly (cf. also \cite{ms15}), then it follows that the $C^0$-limit of a sequence of smooth Legendrian sections of $pr_1: E \to M$ is again Legendrian as long as it is a smooth section. 

Still under the assumption that the conformal factors converge uniformly, Rosen and Zhang \cite{rz18} proved a result analogous to the Humilière-Leclercq-Seyfaddini Theorem, namely, that smooth images of coisotropic submanifolds (i.e. $(TL \cap \xi)^{\perp_{d \alpha}} \subseteq TL$, cf. \cite{hua15}) under homeomorphisms that are $C^0$-limits of contactomorphisms are again coisotropic. Usher \cite{ush20} showed that the conclusion of this statement is still true if the conformal factors are only required to be uniformly bounded from below.

Now we want to examine the question under which conditions smooth $C^0$-limits of Legendrian submanifolds are again Legendrian, even if the limit in not induced by a $C^0$-limit of contactomorphisms. It is well-known that any n-dimensional submanifold of a contact manifold $(M^{2n+1}, \xi)$ can be $C^0$-approximated by Legendrian submanifolds as long as there is no topological obstruction (see \cite{em02}, 16.1.3), but we will show that under certain conditions such approximations must have short Reeb chords.

\begin{thm}\label{c0-limit legendre}
Let $(M^{2n+1}, \xi = \ker{\alpha)}$ be a cooriented contact manifold and $L^n$ a closed manifold. Let $f_i:L \to M$ be a sequence of Legendrian embeddings of $L$ into $M$ that $C^0$-converge to an embedding $f = f_\infty:L \to M$. Assume that there exists $\eps > 0$ such that for all $i \in \N$ there are no Reeb chords of length less than $\eps$ going from $f_i(L)$ to itself. If one of the following conditions is satisfied, then $f$ is a Legendrian embedding.

(a)  The Reeb vector field is nowhere tangent to $f(L)$ and there exist real numbers $a,b \in \R,~a < b$, and a geometrically bounded symplectic manifold $(N, \omega)$ such that there exists a symplectic embedding $i:(M \times [a,b], d(e^s \alpha)) \to (N,\omega)$.

(b) The Reeb vector field is nowhere tangent to $f(L)$ and $M$ is either compact or the contactization\footnote{In fact, we only have to require that the contact form on $M = P \times \R$ is equal to the standard contact form on $P \times \R$ outside of a compact set.}
 $M = P \times \R$ of a Liouville manifold $P$.

(c) $M$ is the contactization $M = P \times \R$ of an exact, geometrically bounded symplectic manifold $P$.\\
\end{thm}

By a Liouville manifold $P$ we mean an open exact symplectic manifold that contains a compact domain $\overline{P} \subset P$ such that the Liouville vector field is transverse to $\partial \overline{P}$, and such that the Liouville flow $\Phi_t$ satisfies $P \setminus \overline{P} = \bigcup_{t>0} \Phi_t(\partial \overline{P})$.\\

\begin{rmk}
(1) Because the question whether the $C^0$-limit is Legendrian does not depend on the contact form, the theorem should be read as, ``If there exists a contact form such that there is a positive uniform lower bound on the length of the Reeb chords of the $f(L_i)$, then the limit is Legendrian". 

(2) It is known that a Liouville manifold is always geometrically bounded. Therefore, $(c)$ immediately implies $(b)$ in the case that $M$ is the contactization of a Liouville manifold. We explicitly stated that part of $(b)$ nonetheless because the proofs of $(b)$ and $(c)$ rely on different results about Legendrian and non-Legendrian submanifolds.

(3) An embedding $i: M \times [a,b] \to N$ as in $(a)$ exists if $(M, \alpha)$ is a boundary component of a compact symplectic manifold with boundary of contact type. 

(4) The fact that non-Legendrian submanifolds can be $C^0$-approximated by Legendrian submanifolds also shows that the closedness condition on $L$ in Theorem \ref{c0-limit lagrange general} cannot be removed. Indeed, a $C^0$-converging sequence of Legendrian submanifolds lifts in the symplectization to a $C^0$-converging (in the weak topology) sequence of cylindrical Lagrangian submanifolds and the limit of the latter sequence is Lagrangian if and only if the limit of the former sequence is Legendrian.\\
\end{rmk}

\begin{rmk}
Now let us consider $C^0$-approximations of paths in $\R^3$ with its standard contact structure $\xi = \ker(dz - y dx)$. 

On the one hand, if the path is induced by the Reeb flow, then it is easy to see that it can be $C^0$-approximated by Legendrians that do not have any Reeb chords. For example, if $L$ is the interval
\begin{equation}
L \coloneqq I = \{(0,0,z) \in \R^3|\ z \in [0,1]\},
\end{equation}
then $L$ can be $C^0$-approximated by Legendrians, whose Lagrangian projection looks like a spiral (Figure \ref{fig:lagrangian projection of legendrian approximating interval}). 

\begin{figure}[h]
\centering
\includegraphics[width=0.3\linewidth]{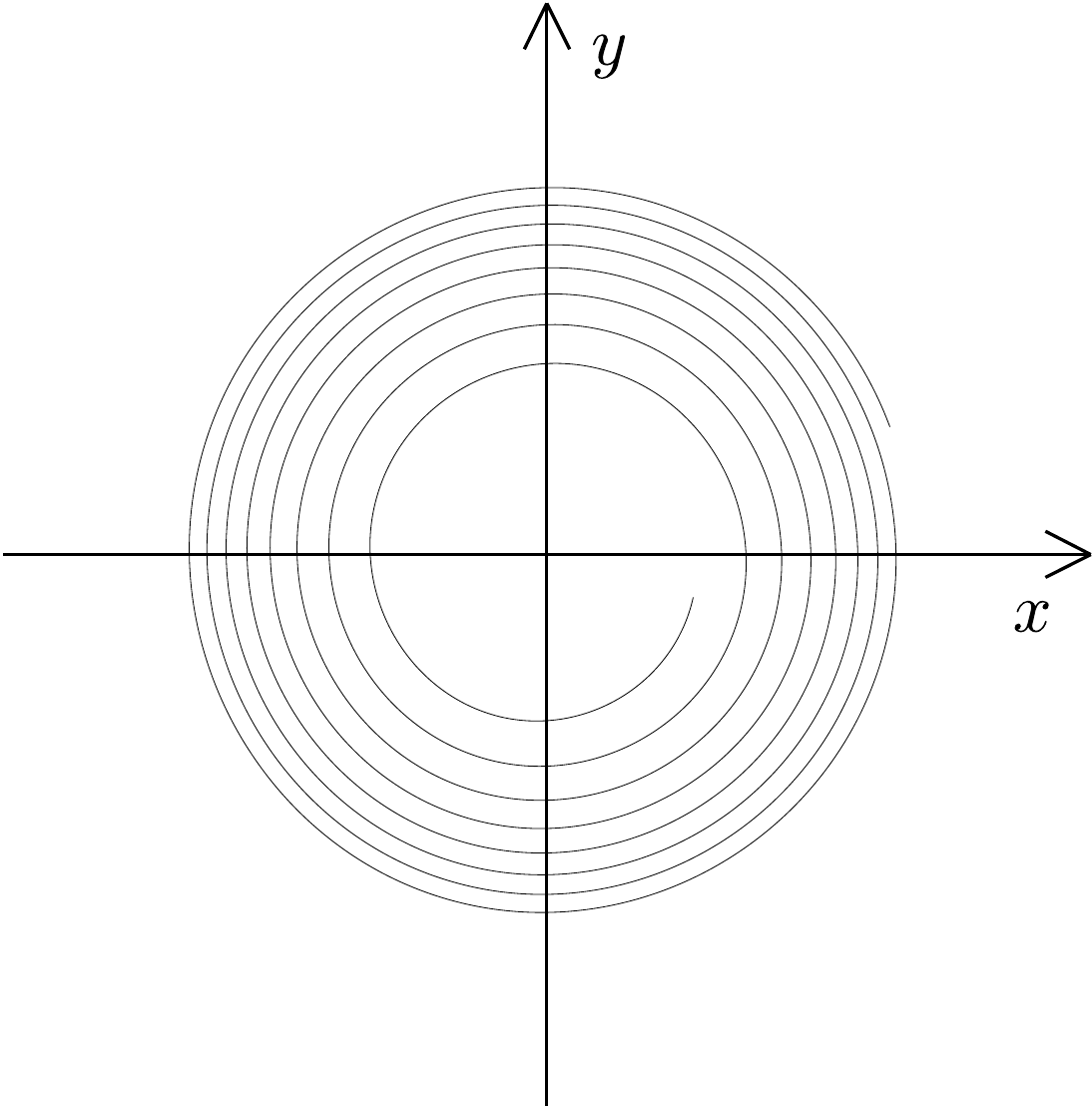}
\caption{Lagrangian projection of a Legendrian submanifold that is $C^0$-approximating the interval.}\label{fig:lagrangian projection of legendrian approximating interval}
\end{figure}

On the other hand, if an embedded path $\gamma: [0,1] \to \R^3$ is not Legendrian and if the Reeb vector field is nowhere collinear to its velocity vector, then it cannot be $C^0$-approximated by a Legendrian path without Reeb chords. In order to see this, let\footnote{We assume that $\gamma$ is defined on the interval $[-1,2]$ instead of $[0,1]$ in order to make it easier to write down the argument below.} $\gamma: [-1,2] \to \R^3$ be such an embedded non-Legendrian path and let $\eta: [-1,2] \to \R^3$ be a Legendrian embedding without Reeb chords that is $\eps$-close to $\gamma$ for some $\eps >0$. Let $\pi: \R^3 \to \R^2$ denote the Lagrangian projection. We write $\pi \gamma$ for $\pi \circ \gamma$ and $\pi \eta$ for $\pi \circ \eta$. $\pi \gamma: [-1,2] \to \R^2$ is an immersed path. Let $\widetilde{\gamma}: [-1,2] \to \R^3$ be the unique Legendrian lift of $\pi \gamma$ to $\R^3$ such that $\widetilde{\gamma}(0) = \gamma(0)$. Since $\gamma$ is not Legendrian, we can assume that, after possibly restricting to a subinterval of $[-1,2]$, $\pi \gamma$ is an embedding and that $\widetilde{\gamma}(1) \neq \gamma(1)$. Let $z$, $\widetilde{z}$ and $z_\eta$ denote the $z$-coordinates of $\gamma$, $\widetilde{\gamma}$ and $\eta$, respectively. Define $C \coloneqq \frac{1}{10}| z(1) - \widetilde{z}(1) | > 0$. Note that $\widetilde{z}$ (and, in fact, the $z$-coordinate of any Legendrian path) satisfies
\begin{equation}\label{eq:z coord of legendrian path}
\widetilde{z}(1) - \widetilde{z}(0) = \int_{\widetilde{\gamma}|_{[0,1]}} y dx.
\end{equation}

For any $\kappa > 0$, let $U_\kappa$ denote the closed $\kappa$-neighbourhood of $\pi \gamma([0,1])$. After possibly decreasing $\eps$, we can assume that there exists a closed ball-shaped neighbourhood $V_\eps$ of $\pi \gamma([0,1])$ that satisfies $U_\eps \subseteq V_\eps \subseteq U_{2\eps}$. We define
\begin{align}
t_- \coloneqq \inf \{t \in [-1,0] | \pi \gamma (s) \in V_\eps ~ \forall s \in [t,0] \}, \\
s_- \coloneqq \inf \{t \in [-1,0] | \pi \eta (s) \in V_\eps ~ \forall s \in [t,0] \},
\end{align}
and similarly we define
\begin{align}
t_+ \coloneqq \sup \{t \in [1,2] | \pi \gamma (s) \in V_\eps ~ \forall s \in [1,t] \}, \\
s_+ \coloneqq \sup \{t \in [1,2] | \pi \eta (s) \in V_\eps ~ \forall s \in [1,t] \}.
\end{align}

Since $\gamma$ and $\eta$ are embedded paths, it follows that $t_-, s_- \to 0$ and $t_+, s_+ \to 1$ as $\eps \to 0$. Now choose $\eps$ so small and $V_\eps$ in such a way that the following conditions are satisfied:

(a) $\eps < C$

(b) $t_-, s_- > -1, \quad t_+,s_+ < 2$,

(c) $\Vert \gamma(s_-) - \gamma(0) \Vert < C$ and $\Vert \gamma(s_+) - \gamma(1) \Vert < C$,

(d) $\Vert \widetilde{\gamma}(t_-) - \widetilde{\gamma}(0) \Vert < C$ and $\Vert \widetilde{\gamma}(t_+) - \widetilde{\gamma}(1 )\Vert < C$,

(e) For any four points $x^0_+,x^0_-,x^1_+,x^1_- \in \partial V_\eps$ with $\Vert x^0_- - x^1_- \Vert < 11\eps$ and $\Vert x^0_+ - x^1_+ \Vert < 11\eps$ and for any two embedded paths $\sigma_0, \sigma_1: [0,1] \to V_\eps$ with $\sigma_{i}(0)= x^{i}_-$ and $\sigma_{i}(1) = x^{i}_+$ for $i \in \{0,1\}$, we have that \begin{equation}
\Big| \int_{\sigma_0} y dx - \int_{\sigma_1} y dx \Big| < C
\end{equation}

(f) $\pi \gamma (\frac{-4 \eps}{\Vert (\pi \gamma)'(0) \Vert}) \not\in U_{3 \eps}$ and $\pi \gamma (1 + \frac{4 \eps}{\Vert (\pi \gamma)'(0) \Vert}) \not\in U_{3 \eps}$.

(g) For all $t \in \left[\frac{-4 \eps}{\Vert (\pi \gamma)'(0) \Vert},0 \right]$ we have that $\Vert \pi \gamma(t) - \pi \gamma(0) \Vert < 5 \eps$, and for all $t \in \left[1, 1 + \frac{4 \eps}{\Vert (\pi \gamma)'(0) \Vert}\right]$ we have that $\Vert \pi \gamma(t) - \pi \gamma(1) \Vert < 5 \eps$.\\

It is clear that (a)-(d) will be satisfied if $\eps$ is sufficiently small. 

To see that (e) can be satisfied, choose $\eps$ so small and choose $V_\eps$ in such a way that for any two points $y^0,y^1 \in \partial V_\eps$ with $\left| y^0 -y^1 \right| < 11 \eps$ there exists an embedded path $\chi:[0,1] \to V_\eps$ with $\chi(0)=y^0$ and $\chi(1)=y^1$ such that $\left|\int_\chi ydx \right| < \frac{C}{10}$. Furthermore, we assume that $\area (V_\eps) < \frac{C}{10}$. Let $x^0_+,x^0_-,x^1_+,x^1_- \in \partial V_\eps$ be four points and  $\sigma_0, \sigma_1: [0,1] \to V_\eps$ be two paths as in (e).
 
By our assumptions, there exist two embedded paths $\chi_-$ and $\chi_+$ with $\chi_-(0) = x_-^0$, $\chi_-(1) = x_-^1$, $\chi_+(0) = x_+^0$ and $\chi_+(1) = x_+^1$ such that $\left|\int_{\chi_\pm} ydx \right| < \frac{C}{10}$. Now let $\lambda_-, \lambda_+,\lambda_0, \lambda_1: [0,1] \to \partial V_\eps$ be four paths that are embeddings when restricted to $(0,1)$ such that $\lambda_-(0) = x_-^0$, $\lambda_-(1) = x_-^1$, $\lambda_+(0) = x_+^0$, $\lambda_+(1) = x_+^1$, $\lambda_0(0) = x_-^0$, $\lambda_0(1) = x_+^0$, $\lambda_1(0) = x_-^1$ and $\lambda_1(1) = x_+^1$. As $d(ydx) = - dx \wedge dy$, it follows from Stokes' Theorem that 
\begin{equation}
\Big| \int_{\chi_-} y dx - \int_{\lambda_-} y dx \Big| \leq \area (V_\eps)< \frac{C}{10}.
\end{equation}

Similarly, it follows that
\begin{equation}
\Big| \int_{\chi_+} y dx - \int_{\lambda_+} y dx \Big| < \frac{C}{10}, \quad \Big| \int_{\sigma_0} y dx - \int_{\lambda_0} y dx \Big| < \frac{C}{10}, \quad \Big| \int_{\sigma_1} y dx - \int_{\lambda_1} y dx \Big| < \frac{C}{10}.
\end{equation}

The absolute value of the winding number of the concatenation $\lambda_0 \ast \lambda_+ \ast \overline{\lambda_1} \ast \overline{\lambda_-}$ is at most four. Here, $\overline{(\cdot)}$ denotes the inversion of paths. Therefore, it follows again from Stokes' Theorem that 
\begin{equation}
\Big| \int_{\lambda_0 \ast \lambda_+ \ast \overline{\lambda_1} \ast \overline{\lambda_-}} y dx \Big| \leq 4~ \area(V) < \frac{4}{10} C.
\end{equation}

Combining the above inequalities one easily concludes that
\begin{equation}
\Big| \int_{\sigma_0} y dx - \int_{\sigma_1} y dx \Big| < C.
\end{equation}

This proves (e).

By looking at the Taylor expansion of $\pi \gamma$ around $0$ and $1$, it can also be seen that (f) and (g) are satisfied if $\eps$ is sufficiently small.

From now on assume that $\eps$ and $V_\eps$ are such that the conditions (a) - (g) are satisfied.

As $\eta$ is $\eps$-close to $\gamma$, (f) implies that $\pi \eta (\frac{-4 \eps}{\Vert (\pi \gamma)'(0) \Vert}) \not\in U_{2 \eps}$ and $\pi \eta (1 + \frac{4 \eps}{\Vert (\pi \gamma)'(0) \Vert}) \not\in U_{2 \eps}$. Since $V_\eps \subseteq U_{2 \eps}$, we can conclude from this observation together with (f) that $s_-, t_- > \frac{-4 \eps}{\Vert (\pi \gamma)'(0) \Vert}$ and $s_+,t_+ < 1 + \frac{4 \eps}{\Vert (\pi \gamma)'(0) \Vert}$. Using (g) and the fact that $\eta$ is $\eps$-close to $\gamma$ it follows that 
\begin{equation}\label{eq:pi eta(s_-)-pi gamma(t_-)}
\Vert \pi \eta (s_-) - \pi \gamma (t_-) \Vert = \Vert \big(\pi \eta (s_-) - \pi \gamma (s_-)\big) + \big(\pi \gamma (s_-) - \pi \gamma(0)\big) + \big(\pi \gamma(0) - \pi \gamma (t_-)\big) \Vert < 11 \eps 
\end{equation}
and similarly also
\begin{equation}\label{eq:pi eta(s_+)-pi gamma(t_+)}
\Vert \pi \eta (s_+) - \pi \gamma (t_+) \Vert < 11 \eps.
\end{equation}

We can see that
\begin{equation}\label{eq:z_eta(0) - z_eta(s_-)}
\left| z_\eta(0) - z_\eta(s_-) \right| = \left| (z_\eta(0) - z(0)) + (z(0) - z(s_-)) + (z(s_-) - z_\eta(s_-)) \right| \leq 3C,
\end{equation}
where in the last inequality we used (a) and (c) together with the assumption that $\eta$ is $\eps$-close to $\gamma$. In the same way we also obtain
\begin{equation}\label{eq:z_eta(s_+) - z_eta(1)}
\left| z_\eta(s_+) - z_\eta(1) \right| \leq 3C.
\end{equation}

We can conclude that
\begin{equation}\label{eq:ztilde diff - z_eta diff}
\begin{gathered}
\left| (\widetilde{z}(1) - \widetilde{z}(0)) - (z_\eta(1) - z_\eta(0)) \right| \\
\overset{\text{(d)},(\ref{eq:z_eta(0) - z_eta(s_-)}), (\ref{eq:z_eta(s_+) - z_eta(1)})}{\leq} \left| (\widetilde{z}(t_+) - \widetilde{z}(t_-)) - (z_\eta(s_+) - z_\eta(s_-)) \right| +  8C \\
\overset{(\ref{eq:z coord of legendrian path})}{=} \Big| \int_{\pi \widetilde{\gamma}|_{[t_-,t_+]}} y dx - \int_{\pi \eta|_{[s_-,s_+]}} y dx \Big| + 8C \overset{\text{(b)},(\ref{eq:pi eta(s_-)-pi gamma(t_-)}),(\ref{eq:pi eta(s_+)-pi gamma(t_+)}),\text{(e)}}{<} 9C,
\end{gathered}
\end{equation}
where in the last step we used the assumptions that $\pi \widetilde{\gamma} = \pi \gamma$ and $\pi \eta$ are embeddings in order to apply (e) (recall that $\eta$ does not have any Reeb chords). 

Now, 
\begin{equation}
\begin{gathered}
|z(1) - z(0) - (z_\eta(1) - z_\eta(0))| \\
\overset{z(0)=\widetilde{z}(0)}{\geq} |z(1) - \widetilde{z}(1)| - |\widetilde{z}(1) - \widetilde{z}(0) - (z_\eta(1) - z_\eta(0))| \overset{\text{def.\ C},~ (\ref{eq:ztilde diff - z_eta diff})}{>} C
\end{gathered}
\end{equation}
leads to a contradiction if $\eps$ is small enough because $\eta$ is $\eps$-close to $\gamma$.

This shows that $\eta$ must have Reeb chords if $\eps$ is sufficiently small. It is also clear that these Reeb chords need to be short because $\eta$ is contained in the $\eps$-neighbourhood of $\gamma$ and the Reeb vector field is nowhere tangent to $\gamma$.

Using Darboux charts, it follows that this statement holds in any $3$-dimensional contact manifold. To the author's knowledge, it is an open question under which conditions it is possible or impossible to $C^0$-approximate open submanifolds $L^n \subseteq (M^{2n+1},\xi = \ker \alpha)$ by Legendrian submanifolds without short Reeb chords in the case $n > 1$.\\
\end{rmk}

\begin{rmk}
If we lower the dimension of $L$ and ask whether the $C^0$-limit of isotropic submanifolds are isotropic, then the answer is no since there is a $C^0$-dense h-principle for subcritical isotropic embeddings into symplectic and contact manifolds (\!\!\cite{em02}, Theorem 12.4.1).\\
\end{rmk}

Another open question is whether Theorem \ref{c0-limit lagrange general} fails if we do not require $M$ to be geometrically bounded, and, similarly, whether the assumptions in Theorem \ref{c0-limit legendre} on $M$ are necessary. Also, one might expect these theorems to hold even for non-compact $L$ if we require the embeddings to be fixed outside some compact subset.\\

\bproofof{Theorem~\ref{c0-limit lagrange general}}
As we can apply the theorem to every connected component of $L$, we can assume that $L$ is connected. 

Recall that for a compactly supported Hamiltonian symplectomorphism $\psi$ the Hofer norm (cf. \cite{hof90}) is defined as
\begin{equation}
\| \psi \| \coloneqq \underset{H} {\inf}\, \|H \|_{osc},
\end{equation}
where the infimum is taken over all time-dependent functions $H_t$ on $M$ whose associated Hamiltonian flow $\phi^H_t$ satisfies $\phi^H_1 = \psi$. Here, $\| H \|_{osc}$ denotes the oscillatory energy of $H$ which is defined as
\begin{equation}
\| H \|_{osc} \coloneqq \int_0^1 \left(\underset{x \in M}{\max}\, H(x,s) - \underset{x \in M}{\min}\, H(x,s)\right) ds.
\end{equation}
The Hofer norm is used to define the displacement energy $e(U)$ of a subset $U \subseteq M$ as
\begin{equation}
e(U) \coloneqq \inf \{\| \psi \| | \psi(U) \cap U = \emptyset\}.
\end{equation}

Now assume that the conclusion of the theorem is false, i.e.\ there exists a sequence of Lagrangian embeddings $f_i:L \to M$ that $C^0$-converge to an embedding $f:L \to M$, but $f$ is not Lagrangian. Let $\iota :S^1 \to T^* S^1$ denote the zero-section. After possibly replacing $L$, $M$, $f_i$ and $f$ by $L \times S^1$, $M \times T^*S^1$, $f_i \times \iota$ and $f \times \iota$, respectively, we can assume that $f(L) \subseteq M$ admits a nowhere-vanishing section of its normal bundle. In order to simplify the notation, we will identify $f(L)$ with $L$ and write $L_i \coloneqq f_i(L)$. By Theorem \ref{main thm symp} there exists a Hamiltonian vector field nowhere tangent to $L$. Hence, for any $\eps > 0$ there is a neighbourhood of $L \subseteq M$ that is displaced by this Hamiltonian isotopy from itself in a time less than $\eps$ by compactness of $L$.

Since the $f_i$ converge uniformly towards $f$, we can find for any $\eps > 0$ a number $N = N(\eps) \in \N$ such that $L_k$ is displaced form itself in a time less than $\eps$ for all $k \geq N$. This implies that the displacement energy of the $L_i$ goes to zero as $i$ increases.

Chekanov proved in \cite{che98} that there is a lower bound on the displacement energy of a closed Lagrangian submanifold in a geometrically bounded symplectic manifold in terms of the minimal area of non-constant pseudoholomorphic spheres in $M$ and non-constant pseudoholomorphic discs in $M$ with boundary on $L$. Let $N \subseteq M$ be a compact tubular neighbourhood of $L$. If $L_i$ is sufficiently $C^0$-close to $L$, then $f_i$ and $f$ are homotopic as maps into $N$. For example, one can explicitly define such a homotopy by moving along the shortest geodesic connecting $f(x)$ and $f_i(x)$ for all $x \in L$. Since $f:L \to N$ is a homotopy equivalence, this implies that $f_i:L \to N$ is a homotopy equivalence as well if $i$ is sufficiently large. Without loss of generality we assume that this is the case for all $f_i$. A non-constant pseudoholomorphic curve with boundary on one of the $L_i$ has positive symplectic area. Hence, it defines a non-trivial class is $H_2(M,L;\R) \cong H_2(M,N;\R)$.

According to Proposition $4.4.1$ in Chapter V in \cite{al94}, there exists a compact neighbourhood $V \subseteq M$ of $N$ such that every pseudoholomorphic curve whose image intersects $N$ lies completely in $V$. Let $U \subseteq M$ be a compact submanifold (possibly with boundary) that contains $V$. Then Lemma \ref{lem:lower bound of area of discs} below shows that the areas of non-constant pseudoholomorphic discs with boundary on one of the $L_i$ are bounded away from zero. Together with Chekanov's energy capacity inequality this implies that the displacement energies of the $L_i$ are uniformly bounded away from zero. This gives the desired contradiction. \eproof\\

\begin{lem}\label{lem:lower bound of area of discs}
Let $N$ be a compact submanifold of a compact manifold $U$ (possibly with boundary). Then there is a constant $C > 0$ such that any disc representing a non-trivial class in $H_2(U,N;\R)$ has area larger than $C$.
\end{lem}

\bproofof{Lemma}
The following proof is an adaptation of the proof of the corresponding lemma in \cite{ls94}. 

By compactness, $H_2(U,N;\Z)$ is finitely generated. Let $\{a_i\}_{i \in I}$ be a basis of the free quotient of $H_2(U,N;\Z)$, where $I$ is a finite index set. Then the $a_i$ also form a real basis of $H_2(U,N;\R)$. Denote by $\{\alpha_i\}_{i \in I}$ the basis of $H^2_{dR}(U,N;\R)$ dual to the $a_i$. Then the homology class of any disc $D$ with boundary on $N$ can be written in the form $D = \sum_{i \in I} n_i a_i$, where $n_i \in \Z$ is the integral of $\alpha_i$ over $D$. Hence, we see that for all $i \in I$,
\begin{equation}
|n_i| = \Big| \int_D \alpha_i \Big| \leq \mathrm{area}(D) \Vert \alpha_i \Vert_{C^0},
\end{equation}
which implies that 
\begin{equation}
\mathrm{area}(D) \geq \mathrm{max}_i \frac{|n_i|}{\Vert \alpha_i \Vert_{C^0}} \geq \mathrm{min}_i \frac{1}{\Vert \alpha_i \Vert_{C^0}},
\end{equation}
if the homology class of $D$ is non-zero (i.e. if not all of the $n_i$ vanish).
\eproof\\

\bproofof{Theorem \ref{c0-limit legendre}} As before we can assume that $L$ is connected.

\textit{Part} $(a)$:
We will reduce Theorem \ref{c0-limit legendre} $(a)$ to Theorem \ref{c0-limit lagrange general} by using the following construction from \cite{moh01}.

First note that since $R_\alpha$ is nowhere tangent to $f_\infty(L)$, we can assume that, after possibly decreasing $\eps > 0$ in the statement of the theorem, there is an $\eps > 0$ such that 
\begin{equation}
\begin{split}
L \times [0, \eps] &\to M \\ 
(x,t) &\mapsto \left( \phi^\alpha_t \circ f_i \right)(x)
\end{split}
\end{equation} 
is an embedding for all $i \in \N \cup \{\infty\}$, where $\phi^\alpha_t$ denotes the Reeb flow in $(M,\alpha)$. Let 
\begin{equation}
(\gamma_1,\gamma_2):S^1 \to [0, \eps] \times [a,b]
\end{equation}
be an embedded loop. Consider the embeddings 
\begin{equation}
\begin{split}
F_i: L \times S^1 & \to (M \times [a,b], d(e^s \alpha)) \\
(x,t) & \mapsto ((\phi^\alpha_{\gamma_1(t)} \circ f_i)(x), \gamma_2(t)).
\end{split}
\end{equation}
It is clear that the $C^0$-convergence of the $f_i$ implies $C^0$-convergence of the $F_i$. Furthermore, a straightforward computation shows that $F_i$ is a Lagrangian embedding if and only if $f_i$ is a Legendrian embedding. Hence, we can apply Theorem \ref{c0-limit lagrange general} to conclude that $f_\infty$ is a Legendrian embedding.\\

The proofs of $(b)$ and $(c)$ are similar to the proof of Theorem \ref{c0-limit lagrange general}. We will use known rigidity results for Legendrian and non-rigidity results for non-Legendrian submanifolds to prove the statement. Again, we identify $L$ with $f(L)$ and write $L_i \coloneqq f_i(L)$.\\

\textit{Part} $(b)$: Assume that $L$ is not Legendrian. After possibly replacing $(M,\alpha)$, $L$, $f_i$ and $f$ by $(M\times T^*S^1,\alpha - p dq)$, $L \times S^1$, $f_i \times \iota$ and $f \times \iota$, respectively, we can assume that there exists a vector field that is nowhere (along $L$) contained in $TL \oplus \< R_\alpha \>$. Here, $\iota:S^1 \to T^*S^1$ denotes the zero section.

Theorem \ref{main thm contact} implies that there exists a contact vector field $X$ that is nowhere contained in $TL \oplus \<R_\alpha\>$. This implies that its flow $\phi_t \coloneqq \phi^X_t$ displaces $L$ for sufficiently small times such that there are no short (compared to the length of the Reeb chords of $L$) Reeb chords between $L$ and $\phi_t(L)$ for any $t>0$ that is sufficiently small. To be more precise, let $\sigma$ denote the minimal length of Reeb chords of $L$. Then, for any $\lambda > 0$ there exists a $\delta > 0$ such that there are no Reeb chords of length smaller than $\sigma - \lambda$ between $L$ and $\phi_t(L)$ for all $0 < t < \delta$. In this case, $\phi_t$ also displaces a neighbourhood of $L \subseteq M$ without short Reeb chords by compactness of $L$. Then for any sufficiently small $t > 0$, there exists an $N \in \N$ such that $\phi_t$ also displaces $L_i$ without short Reeb chords for all $i\geq N$. This shows that for any $\eta > 0$ there is an $N \in \N$ and a function $H:M \to \R$ such that $\Vert H \Vert_{C^1} < \eta$ and the contactomorphism associated to $H$ displaces $L_i$ without short Reeb chords for all $i \geq N$.

For any closed Legendrian submanifold $N \subseteq M$, let $\sigma(\alpha, N)$ denote the minimal length of Reeb chords $\gamma$ of $N$ and of closed Reeb orbits $\gamma$ in $M$ satisfying $[\gamma] = 0 \in \pi_1(M,N)$. Rizell and Sullivan proved that if the $C^1$-norm\footnote{In fact, they only required that the oscillatory energy of $H$ and the conformal factor of the contact flow associated to $H$ are sufficiently small.} of a generic function $H$ on $M$ is small compared to $\sigma(\alpha, N)$, there always exist short (compared to the $C^1$-norm of $H$) Reeb chords between $N$ and $\phi_1^H(N)$ (\!\!\cite{rs16}, Theorem 1.3) if $M$ satisfies the conditions in $(b)$. This gives the desired contradiction because, after possibly approximating $H$, we can assume that it is generic. \\

\textit{Part} $(c)$: For a compactly supported contactomorphism $\psi$ on $(M,\alpha)$ that is isotopic to the identity one can define 
\begin{equation}
\| \psi \|_\alpha \coloneqq \underset{H}{\inf}\, \|H\|,
\end{equation}
where the infimum is taken over all time-dependent functions $H_t$ whose associated contact isotopy $\phi^H_t$ satisfies $\phi^H_1 = \psi$. Here, $\|H \|$ is defined by 
\begin{equation}
\|H\| \coloneqq \int_0^1 \underset{x \in M}{\max}\, H(x,s) ds.
\end{equation}
Shelukhin \cite{she17} proved that this defines a (non-degenerate) norm on the group of compactly supported contactomorphisms isotopic to the identity.

Now assume that $L$ is not Legendrian. Note that for a generic contactomorphism $\phi$, $\phi(L)$ will not intersect $L$ since $\dim(L) = n$ and $\dim(M) = 2n+1$. Let $\phi$ be such a contactomorphism. Then there exists a lower bound $C > 0$ on the length of Reeb chords between $L$ and $\phi(L)$. Theorem 1.9 and Proposition 7.4 in \cite{rz18} together imply that there exist a sequence $\phi_n$ of contactomorphisms isotopic to the identity such that $\lim\limits_{n \to \infty} \Vert \phi_n \Vert_\alpha = 0$ and $\phi_n(L) = \phi(L)$ for all $n \in \N$. By compactness of $L$ we can find for any $n \in N$ and any $\eta > 0$ a neighbourhood $U = U(n,\eta)$ of $L \subseteq M$ such that there are no Reeb chords of length smaller than $C - \eta$ between $U$ and $\phi_n(U)$. After possibly perturbing the $\phi_n$ and choosing a slightly larger $\eta$, we can assume that the $\phi_n$ are generic and still have the above properties (except, of course, $\phi_n(L) = \phi(L)$). 

Since the $L_i$ $C^0$-converge to $L$, we can find for any two positive numbers $\delta, \eta > 0$ some numbers $n,K \in \N$ such that $L_i \subseteq U(n,\eta)$ for all $i \geq K$ and $\Vert \phi_n \Vert_\alpha < \delta$. In particular, there are no Reeb chords of length smaller than $C - \eta$ between $L_i$ and $\phi_n(L_i)$.

This is a contradiction to a result of Rizell and Sullivan \cite{rs18} that states that there have to exist short Reeb chords between $L_i$ and $\phi_n(L_i)$ in the above setting if $\Vert \phi_n \Vert_\alpha$ is sufficiently small. 
\eproof\\

\begin{rmk}
In the proof of $(b)$ we only had to consider Reeb chords $\gamma$ that satisfy $[\gamma] = 0 \in \pi_1(M,f_i(L))$. One could seemingly strengthen the assumption in part $(b)$ of Theorem~\ref{c0-limit legendre} by only requiring that there exists a uniform lower bound on the length of the Reeb chords that satisfy this condition. But it is easy to see that, in fact, compactness of $L$ and the $C^0$-convergence of the $f_i$ imply that there cannot be a sequence of Reeb chords that are non-zero in $\pi_1(M,f_i(L))$ and whose length converges to zero. Indeed, for sufficiently large $i$, $f_i:L \to N$ is a homotopy equivalence between $L$ and a tubular neighbourhood $N$ of $f(L)$ and any sufficiently short Reeb chord of $f_i(L)$ is contained in $N$. Hence, such a Reeb chord is trivial in $\pi_1(M,N) \cong \pi_1(M,f_i(L))$.\\
\end{rmk}

\bibliography{references}
\bibliographystyle{amsalpha}

\end{document}